\documentstyle[12pt,amstex,amssymb,amscd]{article}

\def\Box{\leavevmode\hbox{\vrule\vbox{\hsize=0.6em\hrule
	 \hrule height0.4em depth0.2em width 0pt
	 \hbox{\hskip0.6em plus 0pt minus 1fil}\hrule}\vrule}}
\def\qed{{\unskip\nobreak\hfil\penalty50
	 \hskip2em\hbox{}\nobreak\hfil\Box
	 \parfillskip=0pt \finalhyphendemerits=0 \par
	 \medbreak\noindent\ignorespaces}}

\newcounter{assn}
\renewcommand{\theassn}{(\roman{assn})}
\newenvironment{assertions}%
{\begin{list}{{\normalshape\theassn}\hskip 0.5em}%
{\usecounter{assn}\topsep=0ex \labelwidth=0em \leftmargin=0em
  \labelsep=0em \itemsep=0ex plus 0.2ex \parsep=\itemsep
  \partopsep=0ex}}%
{\end{list}}

\newenvironment{pf*}{{\em Proof\/:}}{\qed}

\newtheorem{cla}{Claim}
\newtheorem{con}[cla]{Conjecture}
\newtheorem{cor}[cla]{Corollary}
\newtheorem{lem}[cla]{Lemma}
\newtheorem{pro}[cla]{Proposition}

\newtheorem{rem}[cla]{Remark}
\newtheorem{rems}[cla]{Remarks}
\newtheorem{thm}[cla]{Theorem}

\newcommand{\be}{\begin{equation}}
\newcommand{\ee}{\end{equation}}

\newcommand{\diag}{\operatorname{diag}}
\newcommand{\Hom}{\operatorname {Hom}}

\newcommand{\Pic}{\operatorname {Pic}}

\newcommand{\Proj}{\operatorname {Proj}}

\newcommand{\SP}{\operatorname{Sp}}
\newcommand{\Spec}{\operatorname {Spec}}
\newcommand{\Sym}{\operatorname{Sym}}

\newcommand{\C}{{\Bbb C}}

\newcommand{\PR}{{\Bbb P}}

\newcommand{\R}{{\Bbb R}}
\newcommand{\Z}{{\Bbb Z}}

\newcommand{\tensor}{\otimes}

\newcommand{\T}{\begin{pmatrix}{\tau_{11}}&{\hdots}&{\tau_{1g}}\\
     {\vdots}&{}&{\vdots}\\{\tau_{g1}}&{\hdots}&{\tau_{gg}}\end{pmatrix}}

\newcommand{\Tanul}{\begin{pmatrix}
     {0}&{\tau_{12}}&{\hdots}&{\hdots}&{\tau_{1,g-1}}\\
     {\tau_{21}}&{0}&&&\vdots\\
     {\vdots}&&&&\vdots\\
     {\vdots}&&&0&{\tau_{g-2,g-1}}\\
     {\tau_{g-1,1}}&{\hdots}&{\hdots}&{\tau_{g-1,g-2}}&{0}\end{pmatrix}}

\newcommand{\e}{{\bold{e}}}
\newcommand{\half}{{\frac 1 2}}
\newcommand{\tr}[1]{{}^t \!{#1}}
\newcommand{\cA}{{\cal A}}
\newcommand{\h}{{\cal H}}
\newcommand{\cL}{{\cal L}}

\newcommand{\cO}{{\cal O}}
\newcommand{\cP}{{\cal P}}
\newcommand{\cR}{{\cal R}}
\newcommand{\cU}{{\cal U}}
\newcommand{\X}{{\Bbb X}}
\newcommand{\Y}{{\Bbb Y}}
\newcommand{\Cast}{{\C^{\ast}}}

\newcommand{\Zmodd}{{\Z/ d\Z}}
\newcommand{\ed}{{\frac{1}{d}}}

\newcommand{\id}{{\bold{1}}}
\newcommand{\bone}{{\bold{1}}}
\newcommand{\od}{\operatorname{d}}

\begin{document}
\pagenumbering{roman}
\title{Very ample linear systems on abelian varieties\\ preliminary version}
\author{O.~Debarre, K.~Hulek, J.~Spandaw}
\thanks{The first author is partially supported by N.S.F. grant DMS 92-03919
and the European Science Project \lq\lq Geometry of
Algebraic Varieties\rq\rq, contract no. SCI-0398-C (A).
The second and third author are supported by the DFG project
\lq\lq Schwerpunktprogramm komplexe Mannigfaltigkeiten\rq\rq.}
\date{}
\maketitle
\tableofcontents

\section{Introduction} \pagenumbering{arabic}

Let $M$ be an ample line bundle on a abelian variety $X$.
(Throughout this paper, we shall work over the complex numbers $\C$.
Hence by abelian variety we always mean complex abelian variety.)
It is known that $M^n$ is very ample for $n\ge3$ (Lefschetz theorem)
and that $M^2$ is very ample if and only if $|M|$ has no fixed
divisor (Ohbuchi theorem). Very little is known for line bundles on $X$
which are not non-trivial powers of ample line bundles.

In this article we consider exclusively the case of line bundles $L$ of type
$(1,\ldots,1,d)$, which are exactly the pullbacks of principal polarizations by
cyclic isogenies of degree $d$.

When $X$ is an abelian surface, such an $L$ is base-point-free if and only if
$d\ge3$
and  $(X,L)$ is not a product \cite[lemma~10.1.2]{BL};
it is very ample if and only if $d\ge5$
and $X$ contains no elliptic curve $E$ such that $L.E\le 2$
\cite[corollary~10.4.2]{BL}.

In higher dimensions, we show that for $(X,L)$ {\em generic} of type
$(1,\ldots,1,d)$ and dimension $g$, the linear system
$|L|$ is base-point-free if and only if
$d\ge g+1$ (see proposition~\ref{75}).
Moreover, the morphism $\phi_L: X \to {\Bbb P}^{d-1}$ that it
defines is birational onto its image if and only if $d\ge g+2$
(see proposition~\ref{76}).

These results are part of a general conjectural picture: for $d>g$, the
morphism $\phi_L$ should be an embedding outside of a set of dimension
$2g+1-d$. In particular, $L$ should be very ample if and only if $d\ge 2g+2$.
(For $g>2$, Barth and Van de Ven have shown that $L$ cannot be very ample for
$d\le 2g+1$ \cite{B,VV}.)

We show that for $(X,L)$ generic of type $(1,\ldots,1,d)$ and dimension $g$,
the line bundle
$L$ is very ample for $d>2^g$, by checking it on a rank-$(g-1)$ degeneration
(see corollary~\ref{171}).
The same result was proved in \cite{BLR}
by a different method for $g=3$ and $d\ge13$.
For $g=3$, this leaves only the case $d=8$ open (the polarization is never
very ample in that case for the degenerations we consider).

As shown by Koll\'ar in \cite{K}, our result in dimension 3 implies the
following
version of a conjecture of Griffiths and Harris: for $d$ odd, $d\ge 9$,
the degree of any curve on a very general hypersurface of degree $6d$ in
${\Bbb P}^4$ is divisible by $d$.

This paper was completed during a visit of the first author at
the University of Hannover. The authors are grateful to the DFG for
financial support which made this visit possible.

\section{Linear systems on abelian varieties}

In this section we focus on the behaviour of morphisms
$\phi_L: A\to {\Bbb P}(H^0(A,L)^{\ast})$, where $A$ is a
{\em generic} abelian variety and $L$ an ample line bundle on $A$,
such that $h^0(A,L)=d$ or,
equivalently, $L^g = g! d$.

\begin{thm}\label{157} Let $A$ be an abelian variety of dimension $g$
and let $\phi: A\to {\Bbb P}^{d-1}$ be a finite morphism.
Then
\begin{assertions}
\item the ramification locus of $\phi$ has dimension at least $2g-d$
\item if $F$ is a closed subset of $A$ such that the restriction
of $\phi$ to $A- F$ is an embedding, then $\dim(F)\ge 2g+1-d$,
except if $(g,d)=(1,3)$ or $(2,5)$.
\end{assertions}
\end{thm}
\begin{pf*} The ramification of $\phi$ is the locus where
$\od\!\phi:TA \to \phi^{\ast} T \PR^{d-1}$ has rank less than $g$.
Since $\phi$ is finite, $\phi^{\ast} T{\Bbb P}^{d-1}$ is ample, hence
so is $\phi^{\ast} T{\Bbb P}^{d-1} \tensor T^{\ast} A$. By
\cite[theorem~1.1]{FL},
this locus is non-empty if $g\ge d-1-(g-1)=d-g$
and has dimension at least $g-(d-g)$. This proves (i).

To prove (ii), we follow ideas of Van de Ven (\cite{VV}).
Assume that $\dim(F) \le 2g-d$. The intersection of $\phi(A)$ with
$s=2g-d+1$ generic hyperplanes is a smooth irreducible $(g-s)$-dimensional
scheme $S$ contained in $\phi(A-F)$.
Note that $S$ sits in ${\Bbb P}^{2(g-s)}$, so that we can do a
Chern class computation as in \cite{VV}.

Set $l=c_1(\cO_S(1))$. Then $c(TS)=(1+l)^{-s}$
and the exact sequence
\[
     0 \to TS \to T{\Bbb P}^{2(g-s)}|_S \to N \to 0
\]
gives
\(
     c(N)=(1+l)^d
\).
It follows (\cite[prop.~3]{VV}) that
\(
     (\deg S)^2 = c_{g-s}(N) =\binom{d}{g-s} l^{g-s}
\).
Since $l^{g-s}=\deg S =\deg A$, we finally get
\(
     \deg A=\binom {d}{g+1}
\).
Since the degree of any ample line bundle on an abelian variety of
dimension $g$ is divisible by $g!$, the conclusion is that $\binom{d}{g+1}$
is divisible by $g!$. We may assume that $g<d\le 2g+1$. We will let
the reader check that this can only happen for $g=1$ or 2.
(For the case $d=2g+1$, see \cite{VV}.)
Cases where
$d=g+1$ are trivially excluded and so is $(g,d)=(2,4)$.

This finishes the proof of (ii).
\end{pf*}

One might expect that there is equality in (i) and~(ii) for a generic
abelian variety, but the answer probably depends on the {\em type}
of the polarization $\phi^{\ast} \cO_{\PR}(1)$. From now on, we will
restrict ourselves to polarizations of type $(1,\ldots,1,d)$.
We shall need the following facts.

Let $L$ be a line bundle of type $(1,\ldots,1,d)$ on an
abelian variety $A$. We define the groups ${\cal G}(L)$ and $K(L)$ as in
\cite[chapter 6]{BL}. Then,
$K(L)$ is isomorphic to $(\Zmodd)^2$, and there is a central extension:
\[
     1 \to \Cast \to {\cal G}(L) \buildrel{\pi}\over{\to} K(L) \to 0.
\]
The group ${\cal G}(L)$ operates on $H^0(A,L)$ \cite[p.~295]{M2}. If
$\tilde\epsilon$ is an
element of ${\cal G}(L)$ of order $d$, it generates a maximal level subgroup
of ${\cal G}(L)$ (in
the sense of \cite[p.~291]{M2}, hence there exists a non-zero section $s$ of
$L$, unique up to
multiplication by a non-zero scalar, such that $\tilde\epsilon\cdot s=s$
\cite[prop.~3]{M2}. In
particular, for $0\le\lambda< d$, there exists a non-zero section
$s_{\lambda}$ such that
$\tilde\epsilon\cdot s_{\lambda}=e^{2i\pi\lambda/d} s_{\lambda}$. Given $A$,
$L$ and $\tilde\epsilon$,
the ordered set $\{ s_0,\ldots ,s_{d-1}\}$ is well-defined up to
multiplication of its elements by
non-zero scalars; it is a basis for $H^0(A,L)$, which we will call a
{\em canonical} basis.

\begin{pro}\label{75} Let $(A,L)$ be a generic polarized abelian variety of
dimension
$g$ and type $(1,\ldots,1,d)$. Then $|L|$ is base-point-free
if and only if $d>g$.
\end{pro}
\begin{pf*} If $d\le g$, then $d$ elements in $|L|$ always intersect,
since $L$ is ample.

Now let $d>g$.
It is enough to exhibit one example.
The construction is the same as in \cite{BLR}.
Let $E_1,\ldots,E_g$
be elliptic curves, let $\epsilon_j$ be a point of order $d$ on $E_j$,
for $j=1,\ldots,g$ and let $\pi: E_1\times\cdots\times E_g\to A$ be the
quotient
by the subgroup generated by $\epsilon_j-\epsilon_k$, for $1\le j<k\le g$.

Let $0_j$ be the origin of $E_j$, let
$L_j={\cal O}_{E_j}(d 0_j)$ and let $M$ be the polarization
$\bigotimes_{j=1}^g {\rm pr}_j^*L_j$ on
$E_1\times\cdots\times E_g$. Pick a lift $\tilde\epsilon_j$ of $\epsilon_j$ of
order $d$ in ${\cal
G}(L_j)$. Then the $\tilde\epsilon_j\tilde\epsilon_1^{-1}$, for $1< j\le g$,
generate a level
subgroup of ${\cal G}(M)$ hence, by \cite[prop.~1]{M2}, there exists a
polarization $L$ on $A$ of type $(1,\ldots,1,d)$
such that $\pi^*L=M$. Moreover,
if $\{ s_{j,\lambda}\} _{0\le\lambda\le d-1}$ is a canonical basis for
$(E_j,L_j,\tilde\epsilon_j)$,
then
$\{ s_{1,\lambda}s_{2,\lambda}\cdots s_{g,\lambda}\} _{0\le\lambda\le d-1}$
is a
basis for $\pi^*H^0(A,L)$....

If all these sections vanish at some point $(e_1,\ldots,e_g)$ and $d>g$,
then at least two different $s_{j \alpha}$, $s_{j \beta}$ must vanish
at $e_j$ for some $j$. But this cannot happen since
$\operatorname{div}(s_{j\beta})$ and
$\operatorname{div}(s_{j\alpha})$ are distinct and are both translates of the
divisor
$\sum_{l=1}^d (l\epsilon_j)$. Hence $L$ is base-point-free.
\end{pf*}

\begin{rems} \normalshape
\begin{assertions}
\item A similar argument shows that for $0<d\le g$, the generic base locus has
dimension exactly $g-d$.
\item By an argument similar to the one used in \cite{M3}, the following
can be shown. Let $\cA$ be a moduli space of polarized abelian varieties
of dimension~$g$ and degree~$g+1$. Then the locus of polarized abelian
varieties for which the corresponding linear system has a base point
is either $\cA$ or a divisor. The proposition implies that it is a
divisor when the type is $(1,\ldots,1,g+1)$. It may happen that this locus
be everything (e.g. for type $(1,2,2)$ (cf. also \cite{NR})).
\end{assertions}
\end{rems}

In view of theorem~\ref{157}, it is tempting to make the following

\begin{con} Let $(A,L)$ be a generic polarized abelian variety
 of dimension $g>2$ and type
$(1,\ldots,1,d)$ with $d>g$ and let $\phi: A\to {\PR}^{d-1}$ be the morphism
associated with $|L|$ (cf. proposition~\ref{75}). Then the ramification of
$\phi$
has dimension $2g-d$ and there exists a closed subset $F$ of $A$
of dimension~$2g+1-d$
such that the restriction of $\phi$ to $A-F$ is an embedding.
\end{con}

It is of course understood that a set of negative dimension is empty.
In particular, the conjecture implies that for $d\ge g+2$, the morphism
$\phi$ should be birational onto its image. This is proven below
(proposition~\ref{76}).
For $d\ge 2g+2$, the line bundle $L$ should be very ample. In \S5, we will
prove that this is the case for $d>2^g$ (see theorem~\ref{71}).

The following proposition shows that, to prove the conjecture for given $d$
and~$g$,
it is enough to exhibity {\em one} polarized abelian variety,
or a suitable degeneration,
for which it holds.

\begin{pro}\label{77} Consider a commutative diagram
\[
\begin{CD}
     X @>f>> Y \\
     @VgVV @VVhV \\
     W @= W,
\end{CD}
\]
where $X$, $Y$, $W$ are analytic (resp. algebraic) varieties and $f$, $g$ and
$h$ are proper.
Let $p$ be a point in $W$, let $X_p$ be the fibre of $f$ over $p$ and assume
that there is a closed (resp. Zariski closed)
subset $F_p$ of $X_p$ such that the restriction of $f$
to $X_p - F_p$ is unramified (an embedding).
Then, for all points $w$ in an open (resp. Zariski open)
neighbourhood of $p$ in $W$, there
exists a closed (resp. Zariski closed)
subset $F_w$ of $X_w$ with $\dim(F_w)\le \dim(F_p)$ such
that the restriction of $f$ to $X_w -F_w$ is unramified (an embedding).
\end{pro}
\begin{pf*} Let $G$ be the support of $\Omega_{X/Y}$.
Since $\Omega_{X/Y} \tensor \cO_{X_p}\cong \Omega_{X_p/Y_p}$ and
$f|_{X_p-F_p}$ is unramified, $G\cap X_p$ is contained in $F_p$.
By semicontinuity of the dimension of the fibres of a morphism
\cite[prop.~3.4, p.~134]{Fi}, this proves the first part of the
proposition, since $f$ is unramified outside of $G$.

Assume now that $f|_{X_p-F_p}$ is an embedding.
Let $Z$ be the union of the components of $X\times_Y X$
other than the diagonal $\Delta_X$, whose image in $W$ contains $p$.
Then, by definition of $\Omega_{X/Y}$, the set $Z\cap\Delta_X$
is contained in $G\times G$; moreover, since $f|_{X_p-F_p}$ is injective,
$Z_p-\Delta_{X_p}$ is contained in $F_p' \times F_p'$, where
$F_p'=f^{-1}(f(F_p))$. Therefore, $Z_p$ is contained in $F_p'\times F_p'$.
Since $F_p'-F_p$ embeds in $f(F_p)$, one has $\dim(F_p')=\dim(F_p)$.
This finishes the proof: one can take for $F_w$ the first projection
of $Z_w$ in $X_w$.
\end{pf*}

\begin{pro}\label{76} Let $(A,L)$ be a generic
polarized abelian variety of dimension~$g$ and
type $(1,\ldots,1,d)$ with $d\ge g+2$. Then the morphism
$\phi: A \to {\PR}^{d-1}$ associated with $|L|$ is birational onto
its image.
\end{pro}
\begin{pf*}
Let ${\cal A}_{g,d,\theta}$ be the moduli space of abelian varieties $A$ of
dimension
$g$, with a polarization $L$ of type $(1,\ldots,1,d)$ and a point
$\tilde\epsilon$ of order $d$ in
${\cal G}(L)$. Let $P=(A,L,\tilde\epsilon)$ be a point
in  ${\cal A}_{g,d,\theta}$, and let $\{ s_0,\ldots ,s_{d-1}\}$
of be a corresponding canonical basis
for $H^0(A,L)$.

We first claim
that for $P$ in a dense open set $U$ of ${\cal A}_{g,d,\theta}$, no $g+1$
distinct
sections in the canonical basis have a common zero. Since
${\cal A}_{g,d,\theta}$ is irreducible
\cite[chapter 8, \S3]{BL}, it is enough to find one point $P$ for which this
property holds. This
follows directly from the proof of proposition 1: keeping the same notation,
the point
$\tilde\epsilon_1$ of ${\cal G}(L_1)$ corresponds to a point $\tilde\epsilon$
of ${\cal G}(L)$ of
order $d$ \cite[prop.~2]{M2}, and the basis for $H^0(A,L)$
given in the proof is canonical for $(A,L,\tilde\epsilon)$; no $g+1$ elements
of this basis vanish
simultaneously.

Now let $(A,L,\tilde\epsilon)$ be an element of $U$ with $A$ {\em simple}
(i.e. such
that $A$ contains no non-trivial abelian subvarieties). The choice of a
canonical basis defines
a morphism
\[
     \phi :A\to \PR ^{d-1}.
\]
Assume that $\phi$ is not birational over its image. Then there exists a
component $D$ of
$A\times_{\PR ^{d-1}}A$, distinct from the diagonal, such that the first
projection $D\to A$ is
surjective. In particular, $\dim (D)\ge g$. Let $m$ be the morphism
\begin{align*}
     m:  D &\to A\\
     (x,y) &\mapsto x-y.
\end{align*}

Let $a$ be a generic point in $m(D)$, let $F$ be an irreducible component of
$m^{-1}(a)$ and let
$G=\operatorname{pr}_2(F)$. Then
\[
     F=\{ (x+a,x) ;\, x\in G\}.
\]
This implies that $\phi(x) =\phi(x+a)$ for all $x\in G$,
hence $L|_G \cong \tau_{-a}^{\ast} L|_G$.
(For any $x\in A$, $\tau_x: A\to A$ denotes translation over $x$.)
Therefore, $a$ lies in the kernel of
\begin{align*}
     A & \to \Pic^0 (\Gamma)\\
     x &\mapsto (\tau_x^{\ast} L \tensor L^{-1})|_\Gamma.
\end{align*}
Since $A$ is
simple, it follows that either $a$ is torsion, in which case $m$ is constant
with image $a$, or else
$m^{-1}(a)$ is finite, in which case $m$ is surjective. In the first case,
all elements of $|L|$ are invariant by translation by $a$, which implies
$a=0$ and contradicts our
choice of $D$. In the second case, there exists $x\in A$ such that
$\phi (x)=\phi(x+\epsilon)$, where
$\epsilon$ is the image of $\tilde\epsilon$ in $K(L)$. Since
\[
     \phi (x+\epsilon)=(s_0(x),\omega s_1(x),\ldots ,\omega^{d-1}s_{d-1}(x)),
\]
where $\omega= e^{2i\pi/d}$, it follows that all of the
$d$ sections in the canonical basis but one vanish at $x$.
When $d\ge g+2$, this contradicts the fact
that $P\in U$. It follows that $\phi$ is birational onto its image.
\end{pf*}

\section{Degeneration of abelian varieties}\label{70}

\subsection{Abelian varieties}\label{50}

We consider the Siegel space of degree $g$, i.e.
\[
     \h_g = \{ \tau \in M(g \times g, \C);\, \tau = \tr \tau,
          \operatorname{Im} \tau >0\}.
\]
Fix an integer $d\ge1$. Then every point
\[
     \tau=\T\in\h_g
\]
defines a $(1,\ldots,1,d)$-polarized abelian variety,
namely
\[
     A_{\tau} = \C^g / L_{\tau},
\]
where $ L_{\tau}$ is the lattice spanned by the rows of the period
matrix
\[
     \Omega_{\tau}=\begin{pmatrix} \tau\\ D \end{pmatrix},\qquad
     D=\diag(1,\ldots,1,d).
\]
We are interested in certain degenerations of these abelian varieties, namely
those which arise if $\tau_{11},\ldots, \tau_{g-1,g-1}$ go to $i \infty$.
We shall first treat the principally polarized case, i.e. $d=1$.
The general case can then be derived easily from this.

We will employ the following notation
\begin{align*}
     \e: \C &\to \Cast\\
          z &\mapsto e^{2\pi i z}.
\end{align*}
Let
\[
     t_{ij}=\e(\tau_{ij}).
\]
Let $z_1,\ldots,z_g$ be the standard coordinates on $\C^g$ and let
\[
     w_i=\e(z_i).
\]
The abelian variety $A_{\tau}$ can be written as a quotient
\[
     A_{\tau}=\Z^g\backslash (\Cast)^g,
\]
where $l=(l_1,\ldots,l_g)\in\Z^g$ operates on $(\Cast)^g$ by
\begin{align*}\label{3}
     l(w_1,\ldots,w_g) &=(w_1',\ldots,w_g')\\
     w_i' &=\prod_{j=1}^g t_{ij}^{l_j}w_i.
\end{align*}
We are interested in what happens as $t_{11},\ldots,t_{g-1,g-1}$ go to zero.

\subsection{Toroidal embedding}

Recall that $\gamma\in \SP(2g,\Z)$ operates on $\h_g$ by
\[
    \gamma = \begin{pmatrix} A & B\\ C & D \end{pmatrix}:
          \tau \mapsto (A \tau + B)(C \tau + D)^{-1}.
\]
(Here $A$, $B$, $C$ and $D$  are $(g \times g)$-matrices.)
The quotient
\(
     \SP(2g,\Z)\backslash \h_g
\)
is the moduli space of principally polarized abelian varieties
of dimension~$g$.
The symplectic group $\SP(2g,\Z)$ contains the lattice subgroup
\[
     \cP = \left\{ \begin{pmatrix} \id & N\\ 0 & \id \end{pmatrix};\,
          N \in \Sym(g, \Z) \right\},
\]
which acts on $\h_g$ by $\tau \mapsto \tau + N$.
We consider the partial quotient
\[
     \cP \backslash \h_g \subset \Sym(g,\Z)\backslash \Sym(g,\C).
\]
Using the coordinates $t_{ij}$, one can make the identification
\[
     \Sym(g,\Z)\backslash \Sym(g,\C) \cong (\Cast)^{\half g(g+1)}.
\]
We use the standard coordinates $T_{ij}$ on
$\Sym(g,\C) \cong \C^{\half g(g+1)}$
and consider the embedding
\[
     \phi_0: (\Cast)^{\half g(g+1)} \to \C^{\half g(g+1)}
\]
given by
\begin{equation}\label{1}
     T_{ij}=\begin{cases}
               \prod_{k=1}^g t_{ki}& \text{if $i=j$}\\
               t_{ij}^{-1} & \text{if $i\neq j$.}
     \end{cases}
\end{equation}
The image of $(\Cast)^{\half g(g+1)}$ under $\phi_0$ is the standard torus
$(\Cast)^{\half g(g+1)}$ in $\C^{\half g(g+1)}$ and the inverse of $\phi_0$
is given by
\begin{equation}
     t_{ij}=\begin{cases}
               \prod_{k=1}^g T_{ki}& \text{if $i=j$}\\
               T_{ij}^{-1} & \text{if $i\neq j$.}
     \end{cases}
\end{equation}
The reason why we are interested in the map $\phi_0$ is that it is
closely related to the so-called principal cone (see \cite[p.~93]{N}),
which plays an essential role in the reduction theory of quadratic forms.
The embedding $\phi_0$ is also a central building block for toroidal
compactifications of $\cA_g$. In fact we have

\begin{lem} The embedding given by $\phi_0$
is the toroidal embedding corresponding
to the principal cone .
\end{lem}
\begin{pf*} Let $n_{ij}\in\Sym(g,\Z)$ be the matrix defined by
\[
     (n_{ij})_{kl}=\begin{cases}
                         1 & \text{if $\{i,j\}=\{k,l\}$}\\
                         0 & \text{otherwise.}
     \end{cases}
\]
The set $\{n_{ij}\}_{1\le i\le j\le g}$ is a basis of $\Sym(g,\Z)$.
We call it the standard basis.
We define another basis $\{n'_{ij}\}_{1\le i\le j\le g}$
by
\[
     n'_{ij}=n_{ii}+n_{jj}-n_{ij}.
\]
This basis defines the principal cone.
In the notation of \cite[p.~5]{O}, we have
\[
     T_{ij}=\e(m'_{ij})\quad\text{and}\quad t_{ij}=\e(m_{ij}),
\]
where $\{ m_{ij}\}$ is the dual of the standard basis
and $\{ m'_{ij}\}$ is the dual of the basis which defines
the principal cone.
The relation between the two bases for $\Sym(g,\Z)$ gives
a relation between the two dual bases and this yields
the required relation between $T_{ij}$ and $t_{ij}$.
\end{pf*}

\subsection{Mumford's construction}

Let
\[
     X_0 = \phi_0(\cP \backslash \h_g) \subset \C^{\half g(g+1)}.
\]
Recall from section \ref{50} that $A_{\tau}$ is the quotient of
$(\Cast)^g$ by the rank-$g$ lattice generated by
\begin{equation}
\begin{aligned}\label{51}
     r_i&=(t_{1i},\ldots,t_{gi})\\
        &=(T_{1i}^{-1},\ldots,\prod_{k=1}^g T_{ki},\ldots,T_{gi}^{-1})
          \qquad\text{($i=1,\ldots,g$).}
\end{aligned}
\end{equation}
This shows that there is a family of abelian varieties $\cA_0 \to X_0$
such that for each point $[\tau]\in X_0$ the fibre $\cA_{[\tau]}$
is isomorphic to $A_{\tau}$.
We now want to add \lq \lq boundary points\rq \rq, i.e. we consider the set
$X$, which is defined to be
the interior of the closure of $X_0$ in $\C^{\half g(g+1)}$.
Mumford's construction enables us to extend the family $\cA_0 \to X_0$
to a family $\cA \to X$ by adding degenerate abelian varieties over the
boundary points.

For this purpose we consider
\[
    A=\C[T_{ij}],\quad I= (T_{ij})\subset A.
\]
Let $K$ be the quotient field of $A$ and consider the torus
\[
     \tilde{G}=\Spec A[w_1,\ldots,w_g,w_1^{-1},\ldots,w_g^{-1}]
          = (\Cast)^g \times \Spec A.
\]
By $\tilde{G}(K)$ we denote the $K$-valued points of $\tilde{G}$.
The character group $\X=\Hom((\Cast)^g,\Cast)$ is spanned by $w_1,\ldots, w_g$.
Finally, we consider the lattice $\Y\subset \tilde{G}(K)$ which is spanned
by the $r_1,\ldots,r_g$ from~(\ref{51}). Following the terminology
of \cite{M}, we shall call $\Y$ the period lattice.

Let
\begin{align*}
     \Phi:\Y &\to \X\\
     r_i &\mapsto w_i.
\end{align*}

\begin{lem}\label{52} The homomorphism $\Phi$ is a polarization in the sense
of Mumford.
\end{lem}
\begin{pf*}
Let $y=\sum y_i r_i$ and $z=\sum z_i r_i\in\Y$. Then the character
$X^{\Phi(y)} =\prod w_i^{y_i}$ given by $\Phi(y)$ satisfies
\begin{equation}\label{2}
\begin{align*}
     X^{\Phi(y)}(z)&=\prod_{i=1}^g \left(
                           \Big( \prod_{j\neq i} T_{ij}^{-z_j} \Big)
                           \Big( \prod_{k=1}^g T_{ki}\Big)^{z_i}
                     \right)^{y_i}\\
                   &= \prod_{i=1}^g T_{ii}^{y_i z_i}\cdot
                      \prod_{i<j} T_{ij}^{(y_i-y_j)(z_i-z_j)}.
\end{align*}
\end{equation}
Hence $X^{\Phi(y)}(z)=X^{\Phi(z)}(y)$ and $X^{\Phi(y)}(y)\in
I$ unless $y=0$.
It follows that $\Phi$ is a polarization in the sense of Mumford.
\end{pf*}

\begin{rem} \normalshape
In fact $\Phi:\Y \to \X$ is an isomorphism, hence it is
a principal polarization in the sense of \cite[chapter II]{HKW}.
\end{rem}

Before we can explain Mumford's construction, we still have to choose
a star $\Sigma\subset\X$.
For $\alpha$, $\beta\in\R^g$, we say that $\alpha\ge\beta$
if $\alpha_i\ge \beta_i$ for $i=1,\ldots,g$.
If furthermore $\alpha\neq \beta$, then we say that $\alpha>\beta$.
Now set
\[
    \Sigma'=\{ \alpha\in\{0,\pm 1\}^g;\, \alpha\ge0\text{ or }\alpha\le0\}.
\]
We identify this set with
\[
     \Sigma =\{X^{\alpha}=\prod w_i^{\alpha_i};\, \alpha\in\Sigma'\}\subset\X.
\]
This is a star in the sense of \cite{M}.

For technical reasons we note

\begin{lem}\label{53} For all
$y\in \Y$ and for all $\alpha \in \Sigma$,
we have $X^{\Phi(y)+\alpha}(y)\in A$.
\end{lem}
\begin{pf*} The calculations in the proof of lemma \ref{52} show that
\[
     X^{\Phi(y)+\alpha_i}(y)=\prod_{i=1}^g T_{ii}^{y_i(y_i+\alpha_i)}\cdot
          \prod_{i<j} T_{ij}^{(y_i-y_j)(y_i-y_j+\alpha_i-\alpha_j)}.
\]
The claim now follows since
$z(z+\beta)\ge0$ if $z\in\Z$ and $\beta\in\{0,\pm 1\}$ and
$\alpha_i$, $\alpha_i-\alpha_j\in \{0,\pm 1\}$
for all $\alpha\in\Sigma$.
\end{pf*}

We are now ready to explain Mumford's construction.
As in \cite{M, HKW} we consider the graded ring
\[
     \cR = \sum_{k=0}^{\infty}
          \{ K[\ldots,X^{\alpha},\ldots]_{\alpha \in \Y} /
          (X^{\alpha + \beta} - X^{\alpha} X^{\beta}, X^0 -1) \} \theta^k,
\]
where $\theta$ is an indeterminate of degree 1 and all other elements
have degree 0. Let $R_{\Phi,\Sigma}$ be the subring of $\cR$ given by
\[
     R_{\Phi,\Sigma} = A[\ldots,
          X^{\Phi(y)+\alpha}(y)X^{2\Phi(y)+\alpha},
          \ldots]_{\alpha\in\Sigma, y\in\Y}
\]
By lemma \ref{53}
\[
     R_{\Phi,\Sigma} \subset A[\ldots,X^{\alpha}\theta,\ldots]_{\alpha\in\Y}.
\]
Let
\[
     \tilde{P}=\Proj R_{\Phi,\Sigma}.
\]
This is a scheme over $\Spec A$. The group $\Y$ acts on $\tilde{P}$
and the desired extension $\cA\to X$ of $\cA_0 \to X_0$ is given by
\[
     \cA = \Y \backslash (\tilde{P}|_X).
\]

The scheme $\tilde{P}$ is covered by the affine open sets
\[
     U_{\alpha,y}=\Spec R_{\alpha,y},
\]
where
\[
     R_{\alpha,y}=A[\ldots,\frac{X^{\Phi(z)+\beta}(z)}{X^{\Phi(y)+\alpha}(y)}
          X^{2\Phi(z-y)+\beta-\alpha},\ldots]_{\beta\in\Sigma, z\in\Y},
\]
as $\alpha$ runs through $\Sigma$ and $y$ runs through $\Y$.
The action of $y\in\Y$ on $\tilde{P}$ identifies $U_{\alpha,0}$ with
$U_{\alpha,y}$, so it suffices to calculate
\begin{gather*}
     U_{\alpha}=U_{\alpha,0}=\Spec R_{\alpha,0}\\
     R_{\alpha}=R_{\alpha,0}=A[\ldots,M_{\beta,z},\ldots]_{\beta\in\Sigma,
     z\in\Y},
\end{gather*}
where
\begin{align*}\label{54}
     M_{\beta,z}&=X^{\Phi(z)+\beta}(z)X^{2\Phi(z)+\beta-\alpha}\\
          &=\prod_{i=1}^g T_{ii}^{z_i(z_i+\beta_i)}\cdot
            \prod_{i<j} T_{ij}^{(z_i-z_j)(z_i-z_j+\beta_i-\beta_j)}\cdot
            \prod_{i=1}^g w_i^{2z_i+\beta_i-\alpha_i}.
\end{align*}

We are not interested in {\em all} degenerations of abelian varieties
arising from this construction, but only in those which correspond to
$\tau_{11},\ldots, \tau_{g-1,g-1} \to i \infty$.
Hence we can fix the entries $\tau_{ij}$ for $i\neq j$ in the matrix $\tau$.
This corresponds to fixing the coordinates $T_{ij}=t_{ij}^{-1}$ for
$i\neq j$ and hence defines an affine subspace
\[
     L=L(\tau_{ij}) \subset \C^{\half g(g+1)}.
\]
We can use $T_{11},\ldots,T_{gg}$ as coordinates on $L$.
For the sake of simplicity, we introduce the notation
\begin{align*}
     \tau_i&=\tau_{ii}\\
     T_i&=T_{ii},\qquad i=1,\ldots,g,
\end{align*}
and consider the ring
\[
     A'=\C[T_1,\ldots,T_g].
\]
By abuse of notation we shall denote the restriction of
$\tilde{P}$ (resp. $U_{\alpha}$) to $L$ also by
$\tilde{P}$ (resp. $U_{\alpha}$). Then
\[
     U_{\alpha}=\Spec R_{\alpha},
\]
where now
\[
     R_{\alpha}=A'[\ldots,M_{\beta,z},\ldots]_{\beta\in\Sigma, z\in\Y},
\]
and
\[
M_{\beta,z}=\prod_{i=1}^g T_{i}^{z_i(z_i+\beta_i)}\cdot
            \prod_{i=1}^g w_i^{2z_i+\beta_i-\alpha_i}.
\]

\begin{pro}\label{55} For $\alpha\in\Sigma$, we have
$R_{\alpha}=A'[X_1,\ldots,X_g,Y_1,\ldots,Y_g]$, where
\[
     X_i=\begin{cases} T_i^{\alpha_i} w_i &\text{if $\alpha\ge0$}\\
                       w_i & \text{if $\alpha\le 0$}
     \end{cases}
\qquad
     Y_i=\begin{cases} w_i^{-1} &\text{if $\alpha\ge0$}\\
                       T_i^{-\alpha_i} w_i^{-1} & \text{if $\alpha\le 0$.}
     \end{cases}
\]
\end{pro}
\begin{pf*} This follows easily from the observation that
$z(z+\beta)-(2z+\beta)\ge -1$ if $z\in\Z$ and $\beta\in\{0,\pm 1\}$.
\end{pf*}

\begin{cor}\label{79} The scheme $\tilde{P}$ is smooth.
\end{cor}
\begin{pf*} It is enough to show that all $U_{\alpha}$ are smooth.
If $\alpha=0$, then $U_{\alpha}= \C^g \times (\Cast)^g$. Now let
$\alpha\neq 0$. We treat the case $\alpha=(1,\ldots,1)$, the other
cases being similar.
Then
\[
     U_{(1,\ldots,1)}=\Spec (\C[T_1,\ldots,T_g,Z_1,\ldots,Z_{2g}] /
     (Z_i Z_{i+g} -T_i)).
\]
Projecting onto $\Spec \C[T_1,\ldots,T_g,Z_1,\ldots,Z_{2g}]$ shows that
$U_{(1,\ldots,1)}\cong \C^{2g}$.
\end{pf*}

Now we consider the set
\[
     V=X \cap L.
\]
One easily checks that $V$ contains the lines where all $T_i$
but one are zero. Let $\tilde{P}_V$ be the restriction of $\tilde{P}$
to $V$.

\begin{pro}
\begin{assertions}
\item The group of periods $\Y$ acts freely and properly discontinuously
on $\tilde{P}_V$.
\item The quotient $\cA_V = \Y \backslash \tilde{P}_V$ is smooth.
Moreover, the family $\cA_V \to V$ is flat. It extends the family
$\cA_0|_{V\cap X_0}$.
In particular, the general fibre is a smooth abelian $g$-fold.
\end{assertions}
\end{pro}
\begin{pf*}
\begin{assertions}
\item This can be done as in \cite[theorem 3.14 (i)]{HKW}.
\item Smoothness follows from (i) and corollary \ref{79}.
The family is flat since $\cA_V$ is smooth of dimension $2g$
and every fibre has dimension $g$. It extends $\cA_0|_{V\cap X_0}$
by construction.
\end{assertions}
\end{pf*}

\subsection{Description of the degenerate abelian varieties}

We now want to describe the fibre of $\cA_V$ over a point
$p=(0,\ldots,0,T_g)$ with $T_g\neq0$.
We shall denote this fibre by $A_p$ and we shall denote the fibres of
$\tilde{P}$, $U_{\alpha}$ and $U_{\alpha,y}$ over this point
by $\tilde{P}_p$ , $(U_{\alpha})_p$ and $(U_{\alpha,y})_p$ respectively.
Recall that
\[
     (U_{\alpha})_p = \Spec (R_{\alpha})_p,
\]
where
\[
     (R_{\alpha})_p = \C [X_1,\ldots,X_g, Y_1,\ldots,Y_g] =
          \C [X_1,\ldots,X_{g-1}, Y_1,\ldots,Y_{g-1}, w_g, w_g^{-1}]
\]
with $X_i$ and $Y_i$ as in proposition \ref{55}.
We first consider $(U_{\alpha})_p$. Clearly $(U_{0})_p = (\Cast)^g$.
In general, $(U_{\alpha})_p$ consists of $2^h$ irreducible components,
where
\[
     h = \# \{ i\, ;\, 1 \le i \le g-1, \alpha_i\neq0 \}.
\]
It is singular for $h\ge1$. Its regular part is the disjoint
union of $2^h$ tori. These tori can be described as follows.
Consider
\[
     w_i = T_i^{-\beta_i},\qquad i=1,\ldots,g,
\]
where $\beta_i\in \{0, -\alpha_i\}$ for $i=1,\ldots,g-1$.
Outside the hyperplanes $\{T_i=0\}$, this defines a section of the torus
bundle
$U_0$. Note that $U_0 = U_{\alpha}$ on $L- \cup_i \{T_i=0\}$.
This section can be extended to a section of $U_{\alpha}$ over $p$,
where it meets exactly one of the $2^h$ tori whose union is the smooth
part of $(U_{\alpha})_p$. This shows also that
$(U_{\alpha})_p \subset (U_{\beta})_p$ if and only if
$\beta\ge\alpha\ge0$ or $\beta\le\alpha\le0$.

The subgroup $\langle r_g \rangle$ of $\Y$ generated by $r_g$
acts on $U_0$
and hence also on $(U_p)_0$. It also acts on $\Spec\C[w_g, w_g^{-1}] = \Cast$
by
\[
     r_g(w_g) = \e(\tau_g) w_g.
\]
The inclusion of rings
\[
     \C[w_g, w_g^{-1}] \subset (R_{\alpha})_p
\]
defines a map
\[
     (U_0)_p \to \Cast,
\]
which is equivariant with respect to the action of
$\langle r_g \rangle$.
In this way we get a semi-abelian variety of rank $g-1$, i.e. an extension
\[
     0 \to (\Cast)^{g-1} \to \Z \backslash (U_0)_p \to E_{\tau_g,1} \to 0,
\]
where $E_{\tau_g,1}$ is the elliptic curve
\[
    E_{\tau_g,1} = \C / (\Z + \Z \tau_g).
\]
The closure of $(U_0)_p$ in $\tilde{P}_p$ has the structure
of a $({\Bbb P}^1)^{g-1}$-bundle over $\Cast$.
Taking the quotient by $\langle r_g \rangle$ gives rise
to a $({\Bbb P}^1)^{g-1}$-bundle over
the elliptic curve $E_{\tau_g,1}$.

We now return our attention to $\tilde{P}_p$. Recall that
\[
     \tilde{P}_p = \cup_{\alpha\in\Sigma,y\in\Y} (U_{\alpha,y})_p.
\]
It follows from our observations above that the regular part
of $\tilde{P}_p$ is the union of countably many tori. These tori
can be labelled in a natural way by elements
$(l_1,\ldots,l_{g-1})\in\Z^{g-1}$: the section given outside the union
of the hyperplanes $\{T_i=0\}$ by
\[
     w_i=T_i^{-l_i}
\]
can be extended over the point $p$ where it meets exactly
one of the tori contained
in $\tilde{P}_p$. We shall label this torus by $(l_1,\ldots,l_{g-1})$.
The element $r_i$ of $\Y$ ($i=1,\ldots,g-1$) then maps the torus
$(l_1,\ldots,l_i,\ldots,l_{g-1})$ to the torus
$(l_1,\ldots,l_i-1,\ldots,l_{g-1})$,
whereas $r_g$ maps each of these tori to itself.

We can now summarize our discussion in

\begin{pro}\label{56} Let $A_p = \Y \backslash \tilde{P}_p$ be the fibre of
$\cA_ V$ over the point $(0,\ldots,0,T_g)$ of~$V$ with $T_g\neq0$.
Then the following holds
\begin{assertions}
\item The regular part $A_p^{\text{reg}}$ of $A_p$ is a semi-abelian
variety of rank~$g-1$. More precisely, there exists an extension
\[
     0 \to (\Cast)^{g-1} \to A_p^{\text{reg}} \to E_{\tau_g,1} \to 0.
\]
\item The normalization of $A_p$ is a  $({\Bbb P}^1)^{g-1}$-bundle over
the elliptic curve
$E_{\tau_g,1}$. The identifications given by the normalization map
are induced by the following identifications on ${(\PR^1)}^{g-1}\times\Cast$
\begin{align*}
     r_i: & ( w_1,\ldots, w_{i-1},\infty, w_i,\ldots, w_g) \mapsto\\
          &\qquad (t_{i1} w_1,\ldots,t_{i,i-1} w_{i-1},0,t_{i,i+1} w_{i+1},
          \ldots,t_{ig} w_g),
\end{align*}
where $i$ runs through $\{1,\ldots,g-1\}$.
\end{assertions}
\end{pro}

\begin{rems}\label{201} \normalshape
\begin{assertions}
\item Let
\[
     a_i = [\tau_{ig}] \in E_{\tau_g,1}\qquad i=1,\ldots,g-1.
\]
Then the action of $r_i$ lies over the translation $x\mapsto x+a_i$ on
the elliptic curve $E_{\tau_g,1}$.
\item The singularities of $A_p$ can be read off from proposition~\ref{55},
but can also be understood in terms of the identifications described
in proposition~\ref{56}~(ii). For every $h\in \{0,\ldots,g-1\}$,
there is a locally closed $(g-h)$-dimensional subset of $A_p$,
where $2^h$ smooth branches meet and where the Zariski tangent space
has dimension $g+h$. The \lq \lq worst\rq \rq\  singularities
of $A_p$ occur
along an elliptic curve isomorphic to  $E_{\tau_g,1}$, where $2^{g-1}$
smooth branches meet.
\end{assertions}
\end{rems}

\subsection{The $(1,\ldots,1,d)$-polarized case.}

We now turn to the case of general $d\ge 1$, i.e. we consider the
period matrix
\[
     \Omega_{\tau}=\begin{pmatrix} \tau\\ D \end{pmatrix},\qquad
     D=\diag(1,\ldots,1,d).
\]
Dividing out by the last $g$ rows of this period matrix gives a torus
$(\Cast)^g$ with coordinates
\[
     w_i=
\begin{cases}
     \e(z_i) & i=1,\ldots,g-1\\
     \e(z_g/d) & i=g.
\end{cases}
\]
Let
\[
     t_{ij}=
     \begin{cases}
     \e(\tau_{ij}) &\text{if $i,j=1,\ldots,g-1$}\\
     \e(\tau_{ig}/d) &\text{if $i=g$ or $j=g$}.
\end{cases}
\]
Then the first $g$ rows of $\Omega_{\tau}$ act on $(\Cast)^g$ by
multiplication by
\[
\begin{cases}
     (t_{i1},\ldots,t_{ig}) &\text{for the $i$-th row, where
     $i=1,\ldots,g-1$}\\
     (t_{1g}^d,\ldots,t_{gg}^d) &\text{for the $g$-th row.}
\end{cases}
\]
Changing the polarization from a principal one to a polarization
of type $(1,\ldots,1,d)$ corresponds to changing the
group of periods $\Y$ to the subgroup
\[
     \Y' = \langle r_1,\ldots, r_{g-1},  r_g^d \rangle.
\]
We shall, therefore, consider the family
\[
     \cA_ V = \Y' \backslash \tilde{P}_ V.
\]
Now the general element is a smooth abelian variety
with a polarization of type $(1, \ldots, 1,d)$. Proposition~\ref{56}
remains unchanged with the one exception that the base curve $E_{\tau_g,1}$
has to be replaced with the elliptic curve
\[
    E_{\tau_g,d} = \C / (\Z d + \Z \tau_g).
\]

\section{Degeneration of the polarization}\label{80}

Here we shall always consider the case of polarizations of type
$(1,\ldots,1,d)$. What we have done in the previous section
was to extend the family $\cA_0|_{V\cap X_0}$
to a family $\cA_V$ over $V$. We would now like to construct
a relative polarization on $\cA_0|_{V\cap X_0}$ which extends to the
family $\cA_V$. Although this can be done, we shall actually do slightly
less: since we are only interested in degenerations belonging to points
$(0,\ldots,0,T_g)$ with $T_g\neq0$, we shall restrict ourselves to
small neighbourhoods of such points.

For each $m\in\R^g$, consider the theta-function
\begin{align*}
     \theta_{m,0}: \h_g\times\C^g &\to \C\\
(\tau,z) &\mapsto \sum_{q\in\Z^g}\e(\half (q+m) \tau \tr (q+m)+(q+m)\tr z),
\end{align*}
Let
\[
     s(\tau)=(-\half \tau_1,\ldots,-\half \tau_{g-1},0)
\]
and
\[
     r=(0,\ldots,0,\ed).
\]
For $k\in\Z$, we can then consider the functions
\begin{align*}
     \theta_k: \h_g\times \C^g &\to \C\\
     (\tau,z) &\mapsto \theta_{kr,0}(\tau,z+s(\tau)).
\end{align*}
Note that this depends only on the class of $k$ in $\Zmodd$.
These functions all have the same automorphy factor and hence
are sections of a line bundle $\cL_\tau$ on $A_\tau$. In
fact, $\cL_\tau$ represents the $(1,\ldots,1,d)$-polarization on $A_\tau$
and the $\theta_k$ define a basis of the space of sections of this line
bundle \cite[p.~75]{I}.
Let
\[
     \tau'=\Tanul
\]
and set
\[
     \tau''=\tr (\tau_{1g},\ldots,\tau_{g-1,g}).
\]

Finally we consider the analogue of the functions $\theta_k$ in
one variable, i.e. the functions
\begin{align*}
     \vartheta_k: \h_1\times\C &\to \C\\
     (\tau,z) &\mapsto \sum_{q\in\Z}\e(\half (q+k/d)^2 \tau +(q+k/d)z),
\end{align*}
which also depend only on $k\in\Zmodd$.

\begin{pro}\label{57} With the notation of \S\ref{70},
the functions $\theta_k$ can be written in the form
\[
     \theta_k(\tau,z)=\sum_{q\in\Z^{g-1}} c_q(\tau')
          \vartheta_k(\tau_g,z_g+q \tau'')
          \prod_{i=1}^{g-1} t_i^{\half q_i(q_i-1)}w_i^{q_i},
\]
with
\[
     c_q(\tau')=\prod_{0<i<j<g} t_{ij}^{q_i q_j}.
\]
\end{pro}
\begin{pf*} This follows from a straightforward computation.
\end{pf*}

\begin{rem} \normalshape
The shift $z\mapsto z+s(\tau)$ was introduced in order
to obtain integer exponents of the variables $t_i$ in the above
description.
\end{rem}

{}From now on, we fix $\tau'$ and $\tau''$. Let $p=(0,\ldots,0,T_g)$
be an element of $V$ with $T_g\neq0$. For small neighbourhoods
$W$ of $p$ in $V$, it follows from proposition~\ref{57}
that we may consider the $\theta_k$ as holomorphic functions on
$W\times (\Cast)^{g-1}\times\C$ (with coordinates
$(T_1,\ldots,T_g;w_1,\ldots,w_{g-1}, z_g)$).

Let $W_0=W\cap X_0$ and let $\cA_{W_0}$, resp. $\cA_W$ be the restriction of
the
family $\cA$ to $W_0$, resp. $W$. Since the automorphy
factors of the functions $\theta_k$ do not depend on
$k$, there exists a line bundle $\cL_0$ on $\cA_{W_0}$ such that the
functions $\theta_k$ are sections of this line bundle.
The line bundle $\cL_0$ defines a relative polarization
on $\cA_{W_0}$. Our aim is to extend the line bundle $\cL_0$
and its sections $\theta_k$ to $\cA_W$.

\begin{pro} The line bundle $\cL_0$ on $\cA_{W_0}$ can be extended to a
line bundle $\cL$ on $\cA_W$. Moreover, the sections of $\cL_0$ defined
by the functions $\theta_k$ can be extended to sections of $\cL$.
\end{pro}
\begin{pf*} Using proposition~\ref{57}, this can be done in the same way
as in \cite[prop.~II.5.13]{HKW} or as in \cite[prop.~4.1.3]{HW}.
Therefore, we shall not give all technical details, but
only an outline of the proof. We first consider the open part $\cU$
of $\tilde{P}_W$ given by the union of the open sets $(U_{0,y})_W$,
where $y\in\Y$. The codimension of the complement of $\cU$ in $\tilde{P}$ is~2.
Each open set $(U_{0,y})_W$ is a trivial torus of rank~$g$ over~$W$.
For $y=0$, we can use coordinates $(T_1,\ldots,T_g;w_1,\ldots,w_{g-1},w_g)$
to identify $(U_0)_W$ with $W\times (\Cast)^g$.
Since the function $\vartheta_k(\tau_g, z_g+q\tau'')$ can be expressed
in terms of the coordinate $w_g=\e(z_g)$,
we can view the functions
$\theta_k$ as functions on $(U_0)_W$. Similarly, using the action of $\Y$,
we can choose coordinates for $(U_{0,y})_W$ for every $y\in\Y$,
such that this open set is also identified with $W\times (\Cast)^g$. In this
way
we can consider the $\theta_k$ also as functions on $(U_{0,y})_W$.
We can think of $\cU$ as a complex manifold which is obtained by glueing the
open
sets $(U_{0,y})_W$.
For every $y\in\Y$, we consider the $\theta_k$ as sections of the trivial
line bundle on $(U_{0,y})_W$. Using the automorphy factors of the functions
$\theta_k$, we can glue these trivial line bundles and obtain a line bundle
$\cL_\cU$ on $\cU$. This can be done in such a way that the action of $\Y$
on $\cU$ lifts to an action of $\Y$ on $\cL_\cU$.
Hence this line bundle descends to a line bundle on $\Y\backslash \cU$.
By construction, the functions $\theta_k$ define sections of this line bundle.
We have now extended the line bundle $\cL_0$ to an open set
of $\cA_W$, whose complement has codimension~2.
To extend the bundle to the whole of $\cA_W$, one can either
use the Remmert-Stein extension theorem (cf. \cite{HKW, HW})
or one can perform a similar construction using all sets
$(U_{\alpha,y})_W$ ($\alpha\in\Sigma$, $y\in\Y$).
\end{pf*}

For future reference we also note

\begin{pro}\label{72} Let $S=\{0,1\}^{g-1}$. Then
\begin{equation*}\label{10}
   \lim_{t_1,\ldots,t_{g-1}\to 0} \theta_k(\tau,z)=
          \sum_{q\in S} c_q(\tau')  \vartheta_k(\tau_g,z_g+q\tau'') w^q,
\end{equation*}
where
\[
     w^q=\prod_{i=1}^{g-1} w_i^{q_i}.
\]
\end{pro}
\begin{pf*} This follows from proposition~\ref{57}.
\end{pf*}
\begin{rem}\label{200} \normalshape
We denote the restriction of the line bundle $\cL$
to $A_p$ by $\cL_p$.
The functions from
proposition~\ref{72} give $d$ sections of $\cL_p$.
\end{rem}

\section{Very ampleness in the case $d>2^g$}

Let $d\ge1$. Given a symmetric matrix
\[
     \tau'=\Tanul,
\]
an element $\tr\tau''=(\tau_{1g},\ldots,\tau_{g-1,g})$ of $\C^{g-1}$
and an element $\tau_g$ of $\h_1$, we have constructed in \S\ref{70}
a degenerate abelian variety $A_p$ of dimension~$g$, whose normalization
is a ${(\PR^1)}^{g-1}$-bundle over the elliptic curve
$E=E_{\tau_g,d}=\C/(\Z d +\Z \tau_g)$. Moreover, there is a commutative
diagram
\[
\begin{matrix}
     E  & \subset & A_p & \stackrel{\phi}{\longrightarrow} & {\PR}^{d-1}\\
     \big\uparrow && \big\uparrow\vcenter{\rlap{$\scriptstyle{\rho}$}} &&
     \big\uparrow\\
     \C\times\{0,\ldots,0\} &\subset & \C\times {(\PR^1)}^{g-1} &
     \stackrel{\Phi}{\longrightarrow} & \C^d,
\end{matrix}
\]
where the map $\Phi$ is defined (with the notation of \S\ref{80}) by
\begin{gather*}
     \Phi = (\phi_0,\ldots,\phi_{d-1})\\
     \phi_k(z_g;w_1,\ldots,w_{g-1})=
          \sum_{q\in S} c_q(\tau')\vartheta_k(\tau_g,z_g+q\tau'')w^q
\end{gather*}
(cf. proposition~\ref{72}). We want to study the rational map $\phi$.
Recall that when $\phi$ is a morphism, $\phi^{\ast}\cO_{\PR^{d-1}}(1)$
is the line bundle $\cL_p$ on $A_p$ defined in
remark~\ref{200} and note that $\phi(E)$ is
a normal elliptic curve of degree $d$ in $\PR^{d-1}$.

For any $z\in\C$, we let $[z]$ be the image of $z$ in $E$ and we set
$a_i=[\tau_{ig}]$ for $i=1,\ldots,g-1$, and $a=\tr (a_1,\ldots,a_{g-1})$.
Recall that $S=\{0,1\}^{g-1}$.

{}From now on, we assume that $\tau''$ is generic. More precisely, it
suffices that
\begin{equation}\label{20}
     \text{{\em the points $a_1,\ldots,a_{g-1}$ of $E$
         are independent over $\Z$}}.
\end{equation}
Then, for any $x\in E$, the subset
\[
     I(x)=\{x+qa;\, q\in S\}
\]
of $E$ has $2^{g-1}$ elements.

The following lemma is a consequence of the Riemann-Roch theorem.

\begin{lem}\label{13} Any set of at most $d-1$ points on $\phi(E)$
is linearly independent.
\end{lem}

\begin{pro}\label{18} Let $\tau''$ be generic.
If $d>2^{g-1}$, then $\phi_0,\ldots,\phi_{d-1}$ have
no common zeroes. In particular, $\phi$ is a morphism
and $\cL_p$ is base-point-free.
\end{pro}
\begin{pf*} For any $Z=(z_g;w_1,\ldots,w_{g-1})$ in $\C\times{(\PR^1)}^{g-1}$,
the vector $\Phi(Z)$ is a linear combination of the vectors
$(\vartheta_k(\tau_g, z_g +q\tau''))_{k\in\Zmodd}$,
whose coefficients are not all zero.
The proposition then follows from lemma~\ref{13}.
\end{pf*}

{}From now on, we assume that $d>2^{g-1}$. Let $z_g\in\C$;
the morphism $\rho$ induces an isomorphism between
$\{z_g\}\times{(\PR^1)}^{g-1}$ and
a closed subscheme of $A_p$\,, which depends only on $x=[z_g]$. We will
denote it
by~$F_x$. Note that $I(x)\subset F_x$.

Since $d>2^{g-1}$, the points of $\phi(I(x))$ are linearly independent.
It follows from proposition~\ref{72}, that the restriction of $\phi$
to $F_x$ is a Segre embedding.

\begin{pro} Let $\tau''$ be generic. If $d>2^g$, then $\phi$  is injective.
\end{pro}
\begin{pf*} From proposition~\ref{56}, we see that the restriction of $\rho$
to $\C\times(\PR^1 -\{\infty\})^{g-1}$ induces a bijection
\[
     B_p\to A_p,
\]
where $B_p$ is an open subset of the normalization of $A_p$,
fibred over $E$ with fibres isomorphic to
$(\PR^1 -\{\infty\})^{g-1}$. For $x\in E$, we let $F_x^0$ be the image in
$A_p$
of the fibre of $x$. It is a subset of $F_x$.

Let $x$, $y$ be two points of $E$.
Since $d>2^g$, the points of $\phi(I(x)\cup I(y))$
are linearly independent. It follows that if $\phi(F_x^0)$ and $\phi(F_y^0)$
meet, then $x\in I(y)$ and $y\in I(x)$.
Condition~(\ref{20}) then implies $x=y$. Hence, for $x\neq y$, the sets
$\phi(F_x^0)$
and $\phi(F_y^0)$ do not meet. Since $\phi|_{F_x}$ is an embedding,
the lemma is proved.
\end{pf*}

\begin{thm}\label{71} Let $\tau''$ be generic.
For $d>2^g$, the morphism $\phi$ is an embedding.
In particular, $\cL_p$ is very ample.
\end{thm}
\begin{pf*}
It remains to prove that the differential is injective
on the Zariski tangent spaces.

After reordering the coordinates, we may assume that we are at a
point
\[
     P=\rho(z_g;0,\ldots,0,v_{h+1},\ldots,v_{g-1}),
\]
where $z_g\in\C$, $v_{h+1},\ldots,v_{g-1}\in\PR^1-\{0,\infty\}$.
The Zariski tangent space of $A_p$ at $P$ has dimension $g+h$
(see remark~\ref{201}).
Moreover, $A_p$ has $2^h$ smooth branches at $P$, which are indexed
by subsets $K$ of $\{1,\ldots,h\}$.
The branch corresponding to $K$ is
\[
     \rho(\{z_g-\tau_K''\}\times(\PR^1)^{g-1}),
\]
where $\tau_K''=\sum_{i\in K} \tau_{ig}$
and in this branch, the point $P_K$ above $P$ is
\[
   \rho(z_g - \tau_K''; v_1',\ldots,v_{g-1}'),
\]
where
\[
     v_l' =
     \begin{cases}
           \infty                        & \text{if $l\in K$}\\
           0                             &
           \text{if $l\in K'=\{1,\ldots,h\}- K$}\\
          \prod_{i\in K} t_{il}^{-1} v_l & \text{for $h<l<g$}
     \end{cases}
\]
(see proposition~\ref{56}~(ii)).

We change the coordinates around $P_K$ by setting
\[
     w_i'' =
     \begin{cases}
           (w_i')^{-1} &\text{if $i\in K$}\\
           w_i'        &\text{otherwise.}
     \end{cases}
\]
In these coordinates, $\Phi$ is given by
\[
     \Phi(z_g;w_1'',\ldots,w_{g-1}'')=
     \Big( \sum_{q\in S} c_q \vartheta_k(\tau_g,z_g +q\tau'')
          (w'')^{\bone -q_K} (w'')^{q_{K'} +q_L}
          \Big)_{k\in\Zmodd},
\]
where
\begin{align*}
     \bone &= (1,\ldots,1)\in S\\
     L &= \{ h+1,\ldots,g-1\}
\end{align*}
and, for $M\subset \{1,\ldots,g-1\}$ and $q\in S$,
we define $q_M\in S$ by
\[
     (q_M)_i =
     \begin{cases}
          q_i &\text{if $i\in M$}\\
          0 &\text{otherwise.}
     \end{cases}
\]

We need to calculate the corresponding Jacobian matrix at $P_K$
(whose new coordinates are
$(z_g-\tau_K''; 0,\ldots,0,v_{h+1}',\ldots,v_{g-1}')$).

\noindent {\em Derivative with respect to $x$.\/}

In the sum, we need only consider indices $q$ such that
\[
     q_i =
     \begin{cases}
          1 &\qquad\text{if $i\in K$}\\
          0 &\qquad\text{if $i\in K'$.}
     \end{cases}
\]
After setting $r=q-\bone_K$, we get
\[
     \Big( \sum\begin{Sb} r\in S\\
          r_1 = \cdots =r_h =0\end{Sb}
          c_{\bone_K +r} \vartheta_k' (\tau_g,z_g+r\tau'') (v')^{r_L}
          \Big)_{k\in\Zmodd}.
\]
Since
\[
     c_{\bone_K +r}= \Big( \prod\begin{Sb} i<j\\ i,j\in K\end{Sb} t_{ij}\Big)
          \Big( \prod\begin{Sb}i\in K\\ j\in L\end{Sb}t_{ij}^{r_{j}}\Big) c_r
\]
and
\[
     (v')^{r_L}=\prod\begin{Sb}i\in K\\ l\in L\end{Sb}t_{il}^{-r_l}v_l^{r_l},
\]
we get a non-zero multiple of the vector
\[
     \Big( \sum\begin{Sb} r\in S\\
          r_1 = \cdots =r_h =0\end{Sb}
          c_{r} \vartheta_k' (\tau_g,z_g+r\tau'') v^{r}
          \Big)_{k\in\Zmodd}.
\]

\noindent {\em Derivative with respect to $w_{\beta}''$, $\beta\in K$.\/}

In the sum, we need only consider indices $q$ such that
\[
     q_i =
\begin{cases}
      1 &\qquad\text{if $i\in K -\{\beta\}$}\\
      0 &\qquad\text{if $i\in K'\cup \{\beta\}$.}
\end{cases}
\]
A similar calculation yields (after setting $r=q-\bone_{K-\{\beta\}}$) a
non-zero multiple of the vector
\[
     \Big( \sum\begin{Sb} r\in S\\
          r_1 = \cdots =r_h =0\end{Sb}
          c_r \vartheta_k(\tau_g, z_g-\tau_{\beta g}+r \tau'')v^r
          \big( \prod_{j\in L} t_{\beta j}^{-r_j} \big)
          \Big)_{k\in\Zmodd}.
\]

\noindent {\em Derivative with respect to $w_{\beta}''$, $\beta\in K'$.\/}

We get a non-zero multiple of
\[
     \Big( \sum\begin{Sb} r\in S\\
          r_1 = \cdots =r_h =0\end{Sb}
          c_r \vartheta_k(\tau_g, z_g+\tau_{\beta g}+r \tau'')v^r
          \big( \prod_{j\in L} t_{\beta j}^{r_j} \big)
          \Big)_{k\in\Zmodd}.
\]

\noindent {\em Derivative with respect to $w_{\beta}''$, $\beta\in L$.\/}

We get a non-zero multiple of
\[
     \Big( \sum\begin{Sb} r\in S\\
          r_1 = \cdots =r_h =0\\
          r_{\beta}=1\end{Sb}
          c_r \vartheta_k(\tau_g, z_g+r \tau'')v^r
          \Big)_{k\in\Zmodd}.
\]

Altogether, letting $K$ vary among all subsets of $\{1,\ldots,h\}$,
we see that the closure of $\od\!\phi (T_P A_p)$ in $\PR^{d-1}$
is spanned by the following $(g+h+1)$ points
\begin{gather*}
     \Big( \sum\begin{Sb} r\in S\\
          r_1 = \cdots =r_h =0\end{Sb}
          c_{r} \vartheta_k (\tau_g,z_g+r\tau'') v^r
          \Big)_{k\in\Zmodd}\\
     \Big( \sum\begin{Sb} r\in S\\
          r_1 = \cdots =r_h =0\end{Sb}
          c_{r} \vartheta_k' (\tau_g,z_g+r\tau'') v^r
          \Big)_{k\in\Zmodd}\\
     \Big( \sum\begin{Sb} r\in S\\
          r_1 = \cdots =r_h =0\end{Sb}
          c_r \vartheta_k(\tau_g, z_g+\epsilon \tau_{\beta g}+r \tau'')v^r
          \big( \prod_{j\in L} t_{\beta j}^{\epsilon q_j} \big)
          \Big)_{k\in\Zmodd}
\end{gather*}
for all $(\beta,\epsilon)\in \{1,\ldots,h\}\times\{-1,1\}$, and
\[
     \Big( \sum\begin{Sb} r\in S\\
          r_1 = \cdots =r_h =0\\
          r_{\beta}=1\end{Sb}
          c_r \vartheta_k(\tau_g, z_g+r \tau'')v^r
          \Big)_{k\in\Zmodd}
\]
for all $\beta\in\{h+1,\ldots,g-1\}$.
Since $d>2^g\ge2^{g-h}(h+1)$, it follows from lemma~\ref{13}
that the $2^{g-h}(h+1)$ vectors
\begin{gather*}
     \big( \vartheta_k(\tau_g, z_g + r\tau'') \big)_{k\in\Zmodd}\\
     \big( \vartheta_k'(\tau_g, z_g + r\tau'') \big)_{k\in\Zmodd}\\
     \big( \vartheta_k(\tau_g, z_g +
          \epsilon \tau_{\beta g} + r\tau'') \big)_{k\in\Zmodd}
\end{gather*}
for $r\in S$ with $r_1=\cdots=r_h=0$, $\epsilon=\pm 1$ and
$\beta\in \{1,\ldots,h\}$, are linearly independent in $\C^d$.

This implies that the $(g+h+1)$ vectors above
are linearly independent in $\C^d$ and proves that the
image of the differential of $\phi$ at $P$ has dimension $g+h$.
The differential of $\phi$ is therefore everywhere injective,
which proves the theorem.
\end{pf*}

\begin{cor}~\label{171} Let $(A,L)$ be a generic polarized abelian variety
of dimension~$g$ and type $(1,\ldots,1,d)$. For $d>2^g$,
the line bundle $L$ is very ample.
\end{cor}

\begin{rem} \normalshape
Similar calculations show that for given $\tau_g$
and $\tau''$ satisfying condition~(4),
and for a generic choice of the matrix $\tau'$,
the morphism $\phi$ is unramified for $d\ge 2^g-g(g-3)/2$
and is an embedding for $d>2^g -g(g-3)/2$.
For $g\ge4$, this improves slightly on the bound in theorem~\ref{71}.
However, $\phi$ is never an embedding for $g=3$ and $d=8$:
for a generic choice of $\tau'$, it is unramified and
identifies (transversally) a finite number of pairs of smooth
points of $A_p$.
\end{rem}

\end{document}